\documentstyle[12pt]{article}
\newcommand{\qed} {\hspace {0.1in} \rule {1.5mm} {3.5mm}}

\newtheorem{tetel}{Proposition}[section]

\newtheorem{theo}{Theorem}

\def\bn{\mbox{\bf N}}

\def\bz{\mbox{\bf Z}}

\def\ni{\vskip 0.1in\noindent}
\def\wk{\widetilde{K}}

\def\ss{\sum(\Gamma,S)}
\def\ssa{\sum(\Gamma,S_1)}
\def\ssb{\sum(\Gamma,S_2)}
\def\z2{\bz_2}
\newcommand{\vi}{\vskip 0.1in \noindent}
\begin{document}
\title{The Euler-characteristic of discrete groups \\
 and Yuzvinskii's entropy
addition formula}
\author{G\'abor Elek}

\maketitle
\vi
\vi
{\em ABSTRACT:} 
We prove that Yuzvinskii's entropy addition formula
can not be extended for the action of certain finitely
generated nonamenable groups.
\vi
\vi
{\em 1991 AMS CLASSIFICATION NUMBERS:} 22D40, 58F03
\vi
\vi
\vi
\vi
\vi
G\'abor Elek
\vi
Mathematical Institute, Hungarian Academy of Sciences
\vi
Budapest 1364, POB 127 Hungary
\vi
e-mail elek@math-inst.hu
\newpage
\section{Introduction}
{\bf Yuzvinskii's entropy addition formula} \cite{LSW},\cite{Y}
\vi
{\em If $\alpha$ is a $\bz^d$-action on the compact group $X$, and
$Y$ is an $\alpha$-invariant normal subgroup of $X$, let $\alpha_Y$
denote the restriction of $\alpha$ to $Y$ and $\alpha_{X|Y}$ be
the induced action on $X/Y$. Then
$$ h(\alpha)=h(\alpha_{X|Y})+h(\alpha_Y)\quad,$$
where $h(\alpha)$ is the topological $\bz^d$-entropy of $\alpha$.}
\vi
\vi
Historically, the notion of topological entropy was introduced
for $\bz$-actions and $\bz^d$-actions.
Later Kieffer \cite{K} extended the definition for amenable groups.
Recently, Friedland introduced a very general definition of entropy
for arbitrary finitely generated group and semigroup actions \cite{F}.
Although the addition formula is not extended yet for finitely
generated amenable groups in general, the work of Ward and Zhang \cite{WZ}
suggests that such extension is possible. We shall show in this
paper that although Friedland's notion of entropy is very nice and
useful, it cannot be used to extend the addition formula for
arbitrary groups and semigroups. In fact, there is no sensible
notion of entropy for which the addition formula can be extended
for free groups, even if we
consider group actions only on compact, metrizable, Abelian groups.
\section{The axioms for the addition formula in symbolic dynamics}
Let $\Gamma$ be a countable group or semigroup and let $S$
be a finite abelian group. Then we consider as usual the
$\sum(\Gamma,S)$ Bernoulli-shift \cite{CP} of $S$-valued functions
on $\Gamma$. Then $\sum(\Gamma,S)$ is a compact, metrizable, Abelian
group with a natural $\Gamma$-action, where the action on the left
is given by
$$\sigma_\gamma(\tau)(\gamma')=\tau(\gamma'\gamma)\quad.$$
A linear subshift is a $\Gamma$-invariant closed subgroup
of $\ss$. Below we give a system of axioms
for the notion of entropy. We call a group or a semigroup {\it SEA}
(symbolically entropy additive) if the following axioms hold:
\vi
(A1)\quad There is a real valued entropy function $h$ on linear subshifts.
\vi
(A2)\quad $h(\ss)=\log_2|S|$ for any finite Abelian group $S$.
\vi
(A3)\quad If $\phi_1\subset\phi_2$ are linear subshifts then
$0\leq h(\phi_1)\leq h(\phi_2)$.
\vi
(A4)\quad If $\phi_1\subset\ssa,\phi_2\subset\ssb$ and
$T:\phi_1\rightarrow \phi_2$ is a $\Gamma$-invariant homeomorphism,
then $h(\phi_1)=h(\phi_2)$.
\vi
(A5)\quad If a  linear subshift $\phi$ contains only finite many
elements, then $h(\phi)=0$.
\vi
(A6)\quad (The Yuzvinskii-axiom) If $T:\phi_1\rightarrow \phi_2$
is a $\Gamma$-invariant surjective continuous map between linear
subshifts, then
$$h(Ker\,\phi_1)+h(\phi_2)=h(\phi_1)\quad.$$
Note that the notions of entropy mentioned in the Introduction
satisfy the first five axioms. So far we know that $\bn$,$\bz$ and
$\bz^d$ are {\it SEA}. Before stating our main result we
prove two easy propositions.
\begin{tetel}\label{egy} If $\Gamma_1$ is a quotient group of $\Gamma$ and
$\Gamma_1$ is not {\it SEA}, then $\Gamma$ is not {SEA} either.
\end{tetel} 
{\bf Proof}: Let $\pi:\Gamma\rightarrow \Gamma_1$
be a surjective group homomorphism. Then for any
finite Abelian group $S$, we have a natural injective map
$$\pi^*:\sum(\Gamma_1,S)\rightarrow \ss$$
defined by $\pi^*(\tau)(\gamma)=\tau(\pi(\gamma))$.
Now let us notice that for any $\gamma\in\Gamma$,
$$\pi^*(\sigma_{\pi(\gamma)}(\tau))=\sigma_\gamma(\pi^*(\tau))\quad.$$
Indeed,
$$\pi^*(\sigma_{\pi(\gamma)}(\tau))(\gamma')=\sigma_{\pi(\gamma)}(\tau)
(\pi(\gamma'))=\tau(\pi(\gamma'\gamma))=\pi^*(\tau)(\gamma'\gamma)=
\sigma_\gamma(\pi^*(\tau))(\gamma')\quad.$$
Consequently, by continuity, $\pi^*$ maps linear
subshifts into linear subshifts. Let $T:\phi_1\rightarrow\phi_2$
be a continuous group homomorphism between linear $\Gamma_1$-subshifts.
Then we can define 
$$\pi^*T:\pi^*\phi_1\rightarrow \pi^*\phi_2\quad,$$
by $$(\pi^*T)(\pi^*\tau)=\pi^*(T(\tau))\quad.$$
It is easy to see that if $T$ is surjective
then $\pi^*T$ is surjective as well. Also,
$\pi^*$ defines an isomorphism between
$Ker\,T$ and $Ker\,(\pi^*T)$. Now suppose that
$\Gamma$ is {\it SEA}. Then we define $h_{\Gamma_1}$
by $h_{\Gamma_1}(\phi)=h_\Gamma(\pi^*(\phi))$. By
our observations above, $h_{\Gamma_1}$ satisfies all the six axioms.
\qed
\begin{tetel}
The unital free semigroup $W_2$ on two generators  is not {\it SEA}.
\end{tetel}
{\bf Proof:}
Suppose that our six axioms hold. The elements of $W_2$ can be identified
with finite words on the two-letter alphabet $\{a,b\}$.
The unit is the empty word $E$, the multiplication is the concatenation.
Let $S=\{0,1\}$, the cyclic group of two elements.
Let $\phi_a\subset\Sigma(W_2,S)$ be defined the following way.
Let $\tau\in \phi_a$ if and only if $\tau(\gamma)=0$, when $\gamma$
starts with the letter $a$.
Then $\phi_a$ is a linear subshift. Similarly, we can define
$\phi_b$. There is an injective, continuous, $W_2$-invariant map
$T:\sum(W_2,S)\rightarrow \phi_a$,
defined by $T(\tau)(\gamma)=\tau(b\gamma), T(\tau)(E)=0$ and
$T(\tau)(a\gamma')=0$ for any $\gamma'\in W_2$.
By our axioms, $h(\phi_a)=1$. Similarly, $h(\phi_b)=1$. Notice that
$\phi_a\cap\phi_b$ has two elements, hence $h(\phi_a\cap\phi_b)=0$. Therefore,
$$1=h(\sum(W_2,S))=h(\phi_a)+h(\phi_b)=2\quad.$$ Thus we got a
contradiction. \qed
\section{Groups of finite type}
First we recall the notion of groups of finite type \cite{Br}.
Let $\Gamma$ be a finitely generated group, then we call it a group
of finite type if there exists a finite simplicial complex $
B\Gamma$, such that its universal cover $E\Gamma$ is contractible.
It is important to remember that the homotopy type of $B\Gamma$
depends only on $\Gamma$, therefore it is meaningful to
 call the Euler-characteristic of $\Gamma$. We shall denote this integer
by $e(\Gamma)$. Free groups, torsion free nilpotent groups and
 the fundamental groups of negatively curved Riemann manifolds
 are all groups of finite type. Note that groups of finite type are
of torsion-free and finitely presented. The free product of two groups
of finite type is again a group of finite type. Now we can state the main
result of our paper.
\begin{theo}
If $\Gamma$ is a group of finite type and has nonzero Euler-characteristic,
then $\Gamma$ is not a \it{SEA} group.
\end{theo}
Before proving our Theorem it is worth to mention some related facts.
By Stalling's Theorem, amenable groups of finite type always have
zero Euler-characteristic. On the other hand free and surface groups
have negative Euler-characteristic. Note that it is not enough to
prove our Theorem for free groups and then apply Proposition \ref{egy}.
There exists groups with Kazhdan's Property (T), which are of finite
type and nonzero Euler-characteristic. These groups have no free factors.
 \section{The Proof of the Theorem}
Let $\widetilde{K}$ be  an infinite, simplicial
complex with a free and simplicial $\Gamma$-action such that the
quotient space
$K=\widetilde{K}/\Gamma$ is a finite simplicial complex.
 Note that $\widetilde{K}$
does not need to be a universal covering, just a simplicial covering of $K$
with deck-transformation group $\Gamma$.
We shall denote by
 $\wk_i$  the set of $i$-dimensional simplices in $\wk$ and by
 $C^i(\wk)$ is the group of $\z2$-valued functions on $\wk_i$, they are
 also known as the  $\z2$-valued $i$-cochains of $\wk$.
\ni Then we have the $\z2$-cohomology complex as usual, where $\z2$
is the cyclic group of two elements.
$$C^0(\wk)\stackrel{d_0}{\rightarrow} C^1(\wk) 
\stackrel{d_1}{\rightarrow} C^2(\wk)\stackrel{d_2}{\rightarrow}\dots
\stackrel{d_{n-1}}{\rightarrow} C^n(\wk)$$
(here  $n$ is the largest dimension of a simplex in $K$).
The definition of the coboundary maps $d_i$ is very easy in this
case. Simply, if $f\in C^i(\wk)$
and $\sigma\in \wk_{i+1}$, then $df(\sigma)=\sum_{\tau_j}f(\tau_j)$,
where the sum is taken over the $i$-dimensional faces of the simplex $\sigma$.
Then we consider the $\z2$-cohomologies as usual.
$H^0({\wk})=Ker\, d_0,\, H^1(\wk)=Ker\, d_1/ Im\, d_0, \,H^2(\wk)=Ker\,d_2/Im\,d_1$
and so on, finally $H^n(\wk)=C^n(\wk)/Im\, d_{n-1}$.
 Now we use the entropy in order
to define some Betti numbers. 
First note that $C^i(\wk)\cong\Sigma(\Gamma,S_i)$, where $S_i=\sum^{n_i}_{j=1}
\z2$. Here $n_i$ is the number of $i$-simplices in the quotient space
$K$. Thus $h(C^i(\wk))=\log_2|S_i|=n_i$.
Notice that all the kernels and images as well all the $C^i(\wk)$'s are linear
subshifts. Thus the entropy betti numbers can be defined as follows :
$b^0_h(\wk)=h(Ker\,d_0)\,,\,b^1_h(\wk)=h(Ker\,d_1)-h(Im\,d_0)\,,\,
\dots,
b^i_h(\wk)= h(Ker\,d_i)-h(Im\,d_{i-1}),\dots
, b^n_h(\wk)=h(C^n(\wk))-h(Im\,d_{n-1})$. One also have the entropy
Euler-characteristic: $e_h(\wk)=b^0_h(\wk)-b^1_h(\wk)
+b^2_h(\wk)-\dots+(-1)^n b^n_h(\wk)$.
\begin{tetel}
$e_h(\wk)=e(K)$ (i.e. the usual topological Euler-characteristic
of the finite simplicial quotient complex $K$).
\end{tetel}
{\bf Proof:} We follow Cohen's proof for the analogous formula
for the $L^2$-Betti numbers \cite{C}.
 He proved that $e(K)=e_{(2)}(\wk)$, where
$e_{(2)}(\wk)$ denotes the $L^2$-Euler characteristic. (this
was proved in the differential setting by Atiyah himself in his classical
Asterisque paper \cite{At}). By our axioms
 we have the following $n+1$ equations.
\ni$ h(C^0(\wk))=h(Im\,d_0)+h(Ker\,d_0)$
\ni$ h(C^1(\wk))=h(Ker\,d_1))+h(Im\,d_1)$
\ni.
\ni.
\ni$ h(C^{n-1}(\wk))=h(Ker\,d_{n-1}))+h(Im\,d_{n-1}))$
\ni$ h(C^n(\wk))=h(C^n(\wk))$

\ni(The last row is not a misprint !)
Now add up the equations above with alternating sign. This gives the
proof of our proposition.\qed\ni

Now let us turn back to the proof of the Theorem. 
Let $\Gamma$ be a  group of finite type satisfying the {\it SEA} property.
 We need to prove that
$e(\Gamma)=0$. Consider the $B\Gamma$ finite simplicial classifying space.
 Since  $\Gamma$ is of finite type, the universal
covering $\wk$ is contractible. Hence, all the $\z2$-cohomologies
of $\wk$ vanish, except the zeroth. Here, we have the constant
functions representing the cohomology space. By our axioms, this linear
subshift has zero entropy. Consequently, all the Betti numbers are
vanishing.
 Thus by our Proposition $e(\Gamma)=0$.\qed

\end{document}